\documentclass[a4paper,12pt]{extarticle}
\usepackage{geometry}
\usepackage[T1]{fontenc}
\usepackage[utf8]{inputenc}
\usepackage[english]{babel}
\usepackage{amsmath}
\usepackage{amsthm}
\usepackage{amssymb}
\usepackage{fancyhdr}
\usepackage{setspace}
\usepackage{graphicx}
\usepackage{colortbl}
\usepackage{tikz}
\usepackage{pgf}
\usepackage{subcaption}
\usepackage{listings}
\usepackage{indentfirst}
\usepackage{graphicx}
\usepackage[
backend=bibtex,
style=numeric,
sorting=none,
maxbibnames=99
]{biblatex}
\bibliography{refs}

\usepackage[colorlinks,citecolor=blue,linkcolor=blue,bookmarks=false,hypertexnames=true, urlcolor=blue]{hyperref}
\usepackage{indentfirst}
\usepackage{mathtools}
\usepackage{booktabs}
\usepackage[flushleft]{threeparttable}
\usepackage{tablefootnote}
\usepackage{pdfpages}
\usepackage{hyperref}

\usepackage{chngcntr}
\counterwithin{table}{section}
\counterwithin{figure}{section}

\graphicspath{{graphics/}}



\addto\captionsrussian{}

\newcommand{\R}{\mathcal{R}}

\colorlet{darkred}{red!80!black}
\colorlet{darkblue}{blue!80!black}
\colorlet{darkgreen}{green!60!black}
\usepackage[english]{babel}
\usepackage{secdot}
\sectiondot{subsection}
\usepackage{color}
\usepackage{verbatim}
\usepackage{float}
\usepackage{mathtools}
\usepackage{todonotes}
\usepackage{amsmath,amsthm,amssymb}
\usepackage{enumerate}
\usepackage{tikz}

\newtheorem{theorem}{Theorem}
\newtheorem{conj}{Observation}

\newtheorem{lemma}[theorem]{Lemma}
\newtheorem{prop}[theorem]{Proposition}

\theoremstyle{definition}
\newtheorem{defi}{Definition} 
\newtheorem{rem}{Remark}

\def\PM#1{*\text{\textup{P}}^{#1}}
\def\NM#1{*\text{\textup{N}}^{#1}}

\def\hx{\hat x}

\title{Experimental Study of the Game Exact Nim$(5,2)$} 
\author{
Vladimir Gurvich\thanks{
National Research University Higher School of Economics (HSE), Moscow, Russia; e-mail:
vgurvich@hse.ru and vladimir.gurvich@gmail.com}
\and
Artem Parfenov\thanks{
National Research University Higher School of Economics (HSE), Moscow, Russia; e-mail: avparfenov@hse.ru and dunno\_o@icloud.com}
\and
Michael Vyalyi\thanks
{National Research University Higher School of Economics, Moscow, Russia;
 e-mail:
vyalyi@gmail.com}
}

\begin{document}
\maketitle

\begin{abstract}
We compare to different extensions of the ancient game of nim: 
Moore's nim$(n, \leq k)$ and exact nim$(n, = k)$. 
Given integers $n$ and $k$  such that  $0 < k \leq n$, 
we consider $n$ piles of stones. 
Two players alternate turns. 
By one move it is allowed to choose  
and reduce any
(i) at most $k$ or (ii) exactly $k$ piles of stones  
in games nim$(n, \leq k)$ and nim$(n, = k)$, respectively.
The player who has to move but cannot is the loser. 
Both games coincide with nim when $k=1$. 
\newline 
Game nim$(n, \leq k)$ was introduced by Moore (1910) 
who characterized its Sprague-Grundy (SG) values 
0 (that is, P-positions)  and  1.
The first open case is SG values 2 for nim$(4, \leq 2)$. 
\newline 
Game nim$(n, = k)$, was introduced in 2018 
% On the Sprague-Grundy function of exact k-Nim,
% Discrete Appl. Math. 239 (2018) 1-14,
% available at http://arxiv.org/abs/1508.04484 ;
% Endre Boros, V. Gurvich, Nhan Bao Ho, Kazuhisa 
% Makino, and Peter Mursi\v{c}.
An explicit formula for its SG function 
was computed for  $2k \geq n$. 
In contrast, case $2k < n$  seems difficult: 
even the P-positions are not known already for nim$(5,=2)$. 
Yet, it seems that 
the P-position of games nim$(n+1,=2)$ and nim$(n+1,\leq 2)$
are closely related. 
(Note that P-positions of the latter are known.)
Here we provide some theoretical and computational 
evidence of such a relation for $n=5$.
\end{abstract}

	{
		\hypersetup{linkcolor=darkblue}
		
	}

\section{Introduction}
An impartial game of two players is defined 
by a directed acyclic graph 
(DAG) whose nodes correspond and directed edges 
 correspond to the positions and moves, respectively. 
 So, the set of possible moves 
 in a fixed position is the same for both players. 
 We assume that the considered DAG is {\em potentially finite}, 
 that is, it may be infinite 
 but the set of positions reachable from a fixed one, 
 by a sequence of successive moves, is finite.
 In particular, each play is finite. 
 This condition itself excludes directed cycles and  
 also implies the existence of the so-called 
 {\em terminal positions} from which there are no moves.

If the initial position is given, moves and positions of the plays starting from this position, form a subgame. We always assume that this subgame is finite. So, one of two players has a winning strategy in this subgame.

        If there is a winning strategy for the player making the move from a position, it is called N-position, otherwise it is called P-position.

        Evaluation of positions can be computed by the following recurrence  (the backward induction). Terminals are P-positions. 
	If there is a move from current position to a P-position, the position is N-position. 
	If all moves from current position lead to N-positions, then the position is P-position.

    The evaluation of a position can be done in time polynomial in the size of the subgame defined by this position. We are interested in games possessing succinct representations of positions. For these games, the algorithmic complexity of the set of P-positions may vary.

    In this paper, we focus on study a game called the exact (5,2) nim. It is a generalization of the game nim and formal definitions are given below. Contrary to  nim, the algorithmic complexity of recognizing P-positions in the exact (5,2) nim is open. Here we present the results of experimental study of the exact (5,2) nim.

    Now we present formal definitions of the games in study.

	\subsection{Nim}
	
	Nim is a game in which there are $n$ piles, each of them contains 
	a positive number of stones. The player must take a positive number of
	stones from any pile. The one who can't do it, lose.

	Let's identify each position with an $n$-dimensional vector $(x_1,\dots, x_n)$ of non-negative integers, $x_i$ is the number of stones in $i$th pile (some piles may be empty during a play). Due to symmetry of the game we assume that coordinates form a non-decreasing sequence. We use this convention for all games considered here.

	\begin{defi}[The Bouton's matrix]
		The Bouton matrix $B(x) = (x_{ij})$ is a binary matrix, in which the $i$-th row is a binary  representation of $x_i$ written starting from the lowest bit.
	\end{defi}

	\begin{theorem}[Criterion of P-positions]\label{classicNIM}
		The position $x$ in nim is \textup{P}-position, if and only if
		\[
		\bigoplus\limits_{i = 1}^{n} x_{ij} = 0\quad\text{for each $j$.}
		\]
	\end{theorem}

	The winning strategy is to make this amount equal to zero on your turn. 
	\cite{Bouton}

	\subsection{The Moore's nim}

The	Moore's nim($n$, $\leq k$) is similar to the nim on $n$ piles, but now a player can take a positive number of stones from at most $k$ piles in one move.

	\begin{defi}[Moore's vector]
		$\displaystyle
			M(x)_j = \sum\limits_{i = 1}^{n} x_{ij} 
		$
	\end{defi}

	\begin{theorem}[Criterion of P-positions \ \cite{Moore}]\label{generalMoore}
		$x$ is a P-position in the Moore's nim ($n$, $\leq k$) if and only if 
		\[
			M(x)_j \equiv 0 \pmod {k + 1}
		\]
  for each $j$.
	\end{theorem}

    \subsection{Exact nim}
	
    	The exact $(n, k)$ nim was defined in \cite{BOROS20181} as a game where two 
    	players take a positive number of stones 
    	from a $k$-element subset of $n$ piles of stones. 
    	A player who cannot take a positive number of stones, 
    	strictly from $k$ piles on his turn, lose. 

        For the games, where $n \leq 2k$, there are criteria in that paper. So, we are interested in case, when $n = 5, k = 2$.

	\section{Connection with Moore's nim}

    We computed a list of positions with the evaluation for each of them up to $(85, 85, 85, 85, 85)$ in the coordinates-wise order. The file with the results can be found by the link \cite{BaseFILE}. Examining this list we have made several observations.  Most of them reveal relations between the exact $(5, 2)$ nim and the Moore's nim $(4, \leq2)$. Namely, a position $x = (x_1, x_2, x_3, x_4, x_5)$, $ x_1 \leq x_2 \leq \dots \leq x_5$, of the exact $(5, 2)$ nim corresponds to a position  $\hat x = (x_1, x_2, x_3, x_4)$ of the Moore's nim $(4, \leq2)$. 
	
	%\subsection{Four types of positions}

Let's classify all positions into four classes: PP, PN, NP, NN. The first letter is the evaluation of position in the exact (5, 2) nim and the second is the evaluation of $\hx$ in the Moore's $(4,\leq2)$ nim.

	Recall that the number of ones in each column of the Bouton's matrix of a P-position in the Moore's nim $(4,\leq2)$ is multiple of 3, i.e. it equals either 0 or 3. We call such columns balanced.
Let $m(x)_j$ be the sum of $\hx_{ij}$, $1\leq i\leq 4$.
	Thus, $i$-th column is balanced if and only if $m(x)_i\in\{0,3\}$.
	
	%Let's call $\PM{}$ positions $(i,j,k,l,m)$ in (5, 2) nim for which
	%$(i,j,k,l)$ - P-position in Moore's nim $(4,\leq2)$.

	%Denote $j$-th column of matrix $x_{*j}$. We need the rows $1 \dots 4$
	%to connect $(5, 2)$ with Moore's $(4, \leq 2)$. %We call it $\hx$,
	%so $\hx_{ij} = x_{ij}$, $1\leq i \leq 4$.

	%\subsection{Observations}

        It turns out that most of the positions in the game are PP or NN — so in most cases the criterion of Moore is useful also for the exact nim. Because of that, we counted the number of positions for which Moore's criterion is mistaken. 
	
	\begin{conj}\label{Moore3}
    		The ratio of $\frac{|PN|}{|PN + PP|}$ positions among all $\PM{}$ only 20\% of all the positions from the list.
	\end{conj}

\begin{figure}[h]
    \includegraphics[scale=0.7]{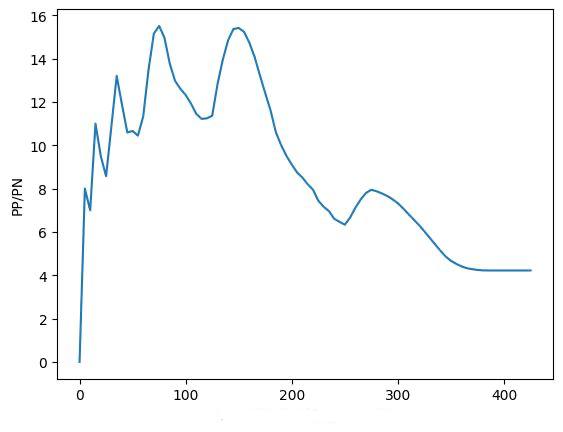}
    \caption{$\frac{|PP|}{|PN|}$ On x-axis, we have the number of stones}\label{pic1}
\end{figure}

    \begin{conj}\label{graphPP/PN}
        $\frac{|PP|}{|PN|}$ isn't monotonous with growth of number of stones, look at Fig \ref{pic1}.
    \end{conj}

On the computed list of the position evaluations, PN-positions reveal an interesting behaviour. 
    %Also, there is an interesting fact about one of the classes of positions.

	%We have a hypothesis about PN-positions.
	Let's write PN-positions in the lexicographic order, starting with the smallest pile. 
Then, if we fix $x_1$, the sequence of  vectors $(x_3-x_2, x_4-x_3, x_5-x_4)$ is periodic. It holds for $x_1\leq 15$ and lengths of periods change in non-monotone way. For larger values of $x_1$ we haven't seen any periods, but it might be caused by the limitations of the computed list.

	\begin{conj}\label{Moore4} 
		$\forall x_1$ vector is periodic.
	\end{conj}

    %For the values $x_1 \leq 15$, we checked on the computer, the hypothesis confirmed.

        According to our notation, in NN and PP positions the evaluations in the exact $(5,2)$ nim and in the Moore's $(4, \leq 2)$  nim coincide, while in NP and PN positions they differ. 

\begin{figure}[h]
    \includegraphics[scale=0.7]{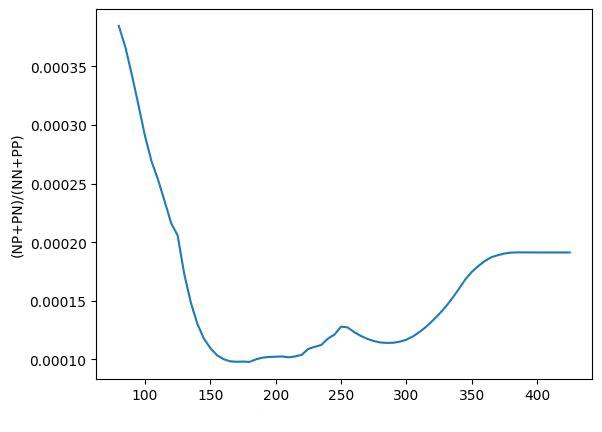}
\caption{$\frac{(|NP| + |PN|)}{(|NN| + |PP|)}$. On x-axis, we have the number of stones.}\label{pic2}
\end{figure}

    \begin{conj}\label{(NP+PN)/(NN+PP)}
        %There are several positions where evaluation in $(5, 2)$ isn't corresponding
		%to evaluation in Moore's $(4, \leq 2)$.
		The plot of the ratio $\frac{(|NP| + |PN|)}{(|NN| + |PP|)}$ at Fig \ref{pic2} shows that in most positions the games are consistent. Also, it shows that the ratio isn't monotonic in the number of stones. 
    \end{conj}

    For the nim and the Moore's nim, there are symmetries of the Bouton's matrix preserving the evaluations of  positions. For example, permutations of columns and insertions/deletions of zero columns do not change the evaluation. We investigated possible symmetries of these sorts for the exact $(5,2)$ nim and did not find any. Below we present a series of counterexamples accompanied by few facts observed. 
    
    %There was an idea that evaluation is connected to transformations of the Bouton's matrix. What if some changes in that matrix such as insertion of the zero-column, permutations of the columns, etc. always change an evaluation, or the opposite — never change it. It turns out that there is no connection, and evaluation might either change or remain the same.
 
	\begin{prop}\label{periodic1}
		There are moves in the exact $(5,2)$ nim such that they result in the permutation of two columns in reduced matrix $\hx{}$ and change the evaluation of positions in the exact (5, 2) nim.	
		For example, 
		\[
		\textcolor{darkblue}{PN} \ (10, 19, 24, 26, 26) \to \textcolor{darkblue}{NN} \ [9, 19, 24, 25, 26] 
		\]
	\end{prop}
	
	\begin{figure}[H]
		\centering
		\begin{tikzpicture}
			\draw[step=1cm,gray,very thin] (0,0) grid (5, 4);
			\draw (0.5,4.5) node[draw, red]{$PN$};
			
			\draw (0.5,0.5) node[draw]{$0$};
			\draw (0.5,1.5) node[draw]{$0$};
			\draw (0.5,2.5) node[draw]{$1$};
			\draw (0.5,3.5) node[draw]{$0$};
			
			\draw (1.5,0.5) node[draw]{$1$};
			\draw (1.5,1.5) node[draw]{$0$};
			\draw (1.5,2.5) node[draw]{$1$};
			\draw (1.5,3.5) node[draw]{$1$};
			
			\draw (2.5,1.5) node[draw]{$0$};
			\draw (2.5,0.5) node[draw]{$0$};
			\draw (2.5,2.5) node[draw]{$0$};
			\draw (2.5,3.5) node[draw]{$0$};
			
			\draw (3.5,1.5) node[draw]{$1$};
			\draw (3.5,0.5) node[draw]{$1$};
			\draw (3.5,2.5) node[draw]{$0$};
			\draw (3.5,3.5) node[draw]{$1$};
			
			\draw (4.5,1.5) node[draw]{$1$};
			\draw (4.5,0.5) node[draw]{$1$};
			\draw (4.5,2.5) node[draw]{$1$};
			\draw (4.5,3.5) node[draw]{$0$};
		\end{tikzpicture}
		\hskip1cm
		\begin{tikzpicture}
			\draw[step=1cm,gray,very thin] (0,0) grid (5, 4);
			\draw (0.5,4.5) node[draw, red]{$NN$};
			
			\draw (0.5,0.5) node[draw]{$1$};
			\draw (0.5,1.5) node[draw]{$0$};
			\draw (0.5,2.5) node[draw]{$1$};
			\draw (0.5,3.5) node[draw]{$1$};
			
			\draw (1.5,0.5) node[draw]{$0$};
			\draw (1.5,1.5) node[draw]{$0$};
			\draw (1.5,2.5) node[draw]{$1$};
			\draw (1.5,3.5) node[draw]{$0$};
			
			\draw (2.5,1.5) node[draw]{$0$};
			\draw (2.5,0.5) node[draw]{$0$};
			\draw (2.5,2.5) node[draw]{$0$};
			\draw (2.5,3.5) node[draw]{$0$};
			
			\draw (3.5,1.5) node[draw]{$1$};
			\draw (3.5,0.5) node[draw]{$1$};
			\draw (3.5,2.5) node[draw]{$0$};
			\draw (3.5,3.5) node[draw]{$1$};
			
			\draw (4.5,1.5) node[draw]{$1$};
			\draw (4.5,0.5) node[draw]{$1$};
			\draw (4.5,2.5) node[draw]{$1$};
			\draw (4.5,3.5) node[draw]{$0$};
			
		\end{tikzpicture}
  \caption{Proposition 3}
	\end{figure}

	\begin{prop}\label{periodic2}
 	There are moves in the exact $(5,2)$ nim such that they result in the permutation of three columns in reduced matrix $\hx{}$ and change the evaluation of positions in the exact (5, 2) nim.	
		For example, 
%		There is a move, which corresponds to the permutation of three columns in $\hx{}$ and changes the evaluation in the exact (5, 2) nim.
	%
		%Example:
		\[
		\textcolor{darkblue}{PN} \ (14, 16, 25, 25, 25) \to \textcolor{darkblue}{NN} \ [7, 8, 25, 25, 25]
		\]
	\end{prop}
	
	\begin{figure}[!h]
		\centering
		\begin{tikzpicture}
			\draw[step=1cm,gray,very thin] (0,0) grid (5, 4);
			\draw (0.5,4.5) node[draw, red]{$PN$};
			
			\draw (0.5,0.5) node[draw]{$1$};
			\draw (0.5,1.5) node[draw]{$1$};
			\draw (0.5,2.5) node[draw]{$0$};
			\draw (0.5,3.5) node[draw]{$0$};
			
			\draw (1.5,0.5) node[draw]{$0$};
			\draw (1.5,1.5) node[draw]{$0$};
			\draw (1.5,2.5) node[draw]{$0$};
			\draw (1.5,3.5) node[draw]{$1$};
			
			\draw (2.5,1.5) node[draw]{$0$};
			\draw (2.5,0.5) node[draw]{$0$};
			\draw (2.5,2.5) node[draw]{$0$};
			\draw (2.5,3.5) node[draw]{$1$};
			
			\draw (3.5,1.5) node[draw]{$1$};
			\draw (3.5,0.5) node[draw]{$1$};
			\draw (3.5,2.5) node[draw]{$0$};
			\draw (3.5,3.5) node[draw]{$1$};
			
			\draw (4.5,1.5) node[draw]{$1$};
			\draw (4.5,0.5) node[draw]{$1$};
			\draw (4.5,2.5) node[draw]{$1$};
			\draw (4.5,3.5) node[draw]{$0$};
		\end{tikzpicture}
		\hskip1cm
		\begin{tikzpicture}
			\draw[step=1cm,gray,very thin] (0,0) grid (5, 4);
			\draw (0.5,4.5) node[draw, red]{$NN$};
			
			\draw (0.5,0.5) node[draw]{$1$};
			\draw (0.5,1.5) node[draw]{$1$};
			\draw (0.5,2.5) node[draw]{$0$};
			\draw (0.5,3.5) node[draw]{$1$};
			
			\draw (1.5,0.5) node[draw]{$0$};
			\draw (1.5,1.5) node[draw]{$0$};
			\draw (1.5,2.5) node[draw]{$0$};
			\draw (1.5,3.5) node[draw]{$1$};
			
			\draw (2.5,1.5) node[draw]{$0$};
			\draw (2.5,0.5) node[draw]{$0$};
			\draw (2.5,2.5) node[draw]{$0$};
			\draw (2.5,3.5) node[draw]{$1$};
			
			\draw (3.5,1.5) node[draw]{$1$};
			\draw (3.5,0.5) node[draw]{$1$};
			\draw (3.5,2.5) node[draw]{$1$};
			\draw (3.5,3.5) node[draw]{$0$};
			
			\draw (4.5,1.5) node[draw]{$1$};
			\draw (4.5,0.5) node[draw]{$1$};
			\draw (4.5,2.5) node[draw]{$0$};
			\draw (4.5,3.5) node[draw]{$0$};
			
		\end{tikzpicture}
  \caption{Proposition 4}
	\end{figure}

 In contrast, we mention the following observation on the computed list.
	
	\begin{conj}\label{periodic3}
		From $NP$ there is no move, which results in the permutation of columns
		of the reduced matrix $\hx{}$. 
	\end{conj}

 Note that in some cases the column premutations and insertions of a zero column do not change the evaluation.
	
	\begin{prop}\label{pereodic4}
		There are  $PN$ positions, which differ by a permutation of two columns in the reduced matrix $\hx{}$. 
For example:
		\[
		\{ (12, 17, 20, 21, 21), (12, 18, 20, 22, 22) \}
		\]
	\end{prop}
	
	\begin{figure}[H]
		\centering
		\begin{tikzpicture}
			\draw[step=1cm,gray,very thin] (0,0) grid (5, 5);
			\draw (0.5,5.5) node[draw, red]{$PN$};
			
			\draw (0.5,0.5) node[draw]{$1$};
			\draw (0.5,1.5) node[draw]{$1$};
			\draw (0.5,2.5) node[draw]{$0$};
			\draw (0.5,3.5) node[draw]{$1$};
			\draw (0.5,4.5) node[draw]{$0$};
			
			\draw (1.5,0.5) node[draw]{$0$};
			\draw (1.5,1.5) node[draw]{$0$};
			\draw (1.5,2.5) node[draw]{$0$};
			\draw (1.5,3.5) node[draw]{$0$};
			\draw (1.5,4.5) node[draw]{$0$};
			
			\draw (2.5,0.5) node[draw]{$1$};
			\draw (2.5,1.5) node[draw]{$1$};
			\draw (2.5,2.5) node[draw]{$1$};
			\draw (2.5,3.5) node[draw]{$0$};
			\draw (2.5,4.5) node[draw]{$1$};
			
			\draw (3.5,0.5) node[draw]{$0$};
			\draw (3.5,1.5) node[draw]{$0$};
			\draw (3.5,2.5) node[draw]{$0$};
			\draw (3.5,3.5) node[draw]{$0$};
			\draw (3.5,4.5) node[draw]{$1$};
			
			\draw (4.5,0.5) node[draw]{$1$};
			\draw (4.5,1.5) node[draw]{$1$};
			\draw (4.5,2.5) node[draw]{$1$};
			\draw (4.5,3.5) node[draw]{$1$};
			\draw (4.5,4.5) node[draw]{$0$};
		\end{tikzpicture}
		\hskip1cm
		\begin{tikzpicture}
			\draw[step=1cm,gray,very thin] (0,0) grid (5, 5);
			\draw (0.5,5.5) node[draw, red]{$PN$};
			
			\draw (0.5,0.5) node[draw]{$0$};
			\draw (0.5,1.5) node[draw]{$0$};
			\draw (0.5,2.5) node[draw]{$0$};
			\draw (0.5,3.5) node[draw]{$0$};
			\draw (0.5,4.5) node[draw]{$0$};
			
			\draw (1.5,0.5) node[draw]{$1$};
			\draw (1.5,1.5) node[draw]{$1$};
			\draw (1.5,2.5) node[draw]{$0$};
			\draw (1.5,3.5) node[draw]{$1$};
			\draw (1.5,4.5) node[draw]{$0$};
			
			\draw (2.5,0.5) node[draw]{$1$};
			\draw (2.5,1.5) node[draw]{$1$};
			\draw (2.5,2.5) node[draw]{$1$};
			\draw (2.5,3.5) node[draw]{$0$};
			\draw (2.5,4.5) node[draw]{$1$};
			
			\draw (3.5,0.5) node[draw]{$0$};
			\draw (3.5,1.5) node[draw]{$0$};
			\draw (3.5,2.5) node[draw]{$0$};
			\draw (3.5,3.5) node[draw]{$0$};
			\draw (3.5,4.5) node[draw]{$1$};
			
			\draw (4.5,0.5) node[draw]{$1$};
			\draw (4.5,1.5) node[draw]{$1$};
			\draw (4.5,2.5) node[draw]{$1$};
			\draw (4.5,3.5) node[draw]{$1$};
			\draw (4.5,4.5) node[draw]{$0$};
		\end{tikzpicture}
  \caption{Proposition 5}
	\end{figure}
	
	\begin{prop}\label{pereodic5}
		There are two $PN$ positions, which differ by insertion of a zero column in the reduced matrix $\hx{}$.
		For example:
		\[
		\{ (6, 9, 10, 11, 11), (12, 17, 20, 21, 21)\}
		\]
	\end{prop}
	
	\begin{figure}[H]
		\centering
		\begin{tikzpicture}
			\draw[step=1cm,gray,very thin] (0,0) grid (5, 5);
			\draw (0.5,5.5) node[draw, red]{$PN$};
			
			\draw (0.5,0.5) node[draw]{$1$};
			\draw (0.5,1.5) node[draw]{$1$};
			\draw (0.5,2.5) node[draw]{$0$};
			\draw (0.5,3.5) node[draw]{$1$};
			\draw (0.5,4.5) node[draw]{$0$};
			
			\draw (1.5,0.5) node[draw]{$1$};
			\draw (1.5,1.5) node[draw]{$1$};
			\draw (1.5,2.5) node[draw]{$1$};
			\draw (1.5,3.5) node[draw]{$0$};
			\draw (1.5,4.5) node[draw]{$1$};
			
			\draw (2.5,0.5) node[draw]{$0$};
			\draw (2.5,1.5) node[draw]{$0$};
			\draw (2.5,2.5) node[draw]{$0$};
			\draw (2.5,3.5) node[draw]{$0$};
			\draw (2.5,4.5) node[draw]{$1$};
			
			\draw (3.5,0.5) node[draw]{$1$};
			\draw (3.5,1.5) node[draw]{$1$};
			\draw (3.5,2.5) node[draw]{$1$};
			\draw (3.5,3.5) node[draw]{$1$};
			\draw (3.5,4.5) node[draw]{$0$};
			
		\end{tikzpicture}
		\hskip1cm
		\begin{tikzpicture}
			\draw[step=1cm,gray,very thin] (0,0) grid (5, 5);
			\draw (0.5,5.5) node[draw, red]{$PN$};
			
			\draw (0.5,0.5) node[draw]{$1$};
			\draw (0.5,1.5) node[draw]{$1$};
			\draw (0.5,2.5) node[draw]{$0$};
			\draw (0.5,3.5) node[draw]{$1$};
			\draw (0.5,4.5) node[draw]{$0$};
			
			\draw (1.5,0.5) node[draw]{$0$};
			\draw (1.5,1.5) node[draw]{$0$};
			\draw (1.5,2.5) node[draw]{$0$};
			\draw (1.5,3.5) node[draw]{$0$};
			\draw (1.5,4.5) node[draw]{$0$};
			
			\draw (2.5,0.5) node[draw]{$1$};
			\draw (2.5,1.5) node[draw]{$1$};
			\draw (2.5,2.5) node[draw]{$1$};
			\draw (2.5,3.5) node[draw]{$0$};
			\draw (2.5,4.5) node[draw]{$1$};
			
			\draw (3.5,0.5) node[draw]{$0$};
			\draw (3.5,1.5) node[draw]{$0$};
			\draw (3.5,2.5) node[draw]{$0$};
			\draw (3.5,3.5) node[draw]{$0$};
			\draw (3.5,4.5) node[draw]{$1$};
			
			\draw (4.5,0.5) node[draw]{$1$};
			\draw (4.5,1.5) node[draw]{$1$};
			\draw (4.5,2.5) node[draw]{$1$};
			\draw (4.5,3.5) node[draw]{$1$};
			\draw (4.5,4.5) node[draw]{$0$};
		\end{tikzpicture}
  \caption{Proposition 6}
	\end{figure}
	
	\begin{prop}\label{pereodic6}
		The permutation of columns in $\hx{}$ might change the evaluation in the exact (5, 2) nim.

		Example:
		
		\[
		\textcolor{darkblue}{PN} \ (12, 17, 20, 21, 21),  \ \textcolor{darkblue}{NN} \ (10, 17, 18, 19, 30) 
		\]
	\end{prop}
	
	\begin{figure}[H]
		\centering
		\begin{tikzpicture}
			\draw[step=1cm,gray,very thin] (0,0) grid (5, 5);
			\draw (0.5,5.5) node[draw, red]{$PN$};
			
			\draw (0.5,0.5) node[draw]{$1$};
			\draw (0.5,1.5) node[draw]{$1$};
			\draw (0.5,2.5) node[draw]{$0$};
			\draw (0.5,3.5) node[draw]{$1$};
			\draw (0.5,4.5) node[draw]{$0$};
			
			\draw (1.5,0.5) node[draw]{$0$};
			\draw (1.5,1.5) node[draw]{$0$};
			\draw (1.5,2.5) node[draw]{$0$};
			\draw (1.5,3.5) node[draw]{$0$};
			\draw (1.5,4.5) node[draw]{$0$};
			
			\draw (2.5,0.5) node[draw]{$1$};
			\draw (2.5,1.5) node[draw]{$1$};
			\draw (2.5,2.5) node[draw]{$1$};
			\draw (2.5,3.5) node[draw]{$0$};
			\draw (2.5,4.5) node[draw]{$1$};
			
			\draw (3.5,0.5) node[draw]{$0$};
			\draw (3.5,1.5) node[draw]{$0$};
			\draw (3.5,2.5) node[draw]{$0$};
			\draw (3.5,3.5) node[draw]{$0$};
			\draw (3.5,4.5) node[draw]{$1$};
			
			\draw (4.5,0.5) node[draw]{$1$};
			\draw (4.5,1.5) node[draw]{$1$};
			\draw (4.5,2.5) node[draw]{$1$};
			\draw (4.5,3.5) node[draw]{$1$};
			\draw (4.5,4.5) node[draw]{$0$};
			
		\end{tikzpicture}
		\hskip1cm
		\begin{tikzpicture}
			\draw[step=1cm,gray,very thin] (0,0) grid (5, 5);
			\draw (0.5,5.5) node[draw, red]{$NN$};
			
			\draw (0.5,0.5) node[draw]{$0$};
			\draw (0.5,1.5) node[draw]{$1$};
			\draw (0.5,2.5) node[draw]{$0$};
			\draw (0.5,3.5) node[draw]{$1$};
			\draw (0.5,4.5) node[draw]{$0$};
			
			\draw (1.5,0.5) node[draw]{$1$};
			\draw (1.5,1.5) node[draw]{$1$};
			\draw (1.5,2.5) node[draw]{$1$};
			\draw (1.5,3.5) node[draw]{$0$};
			\draw (1.5,4.5) node[draw]{$1$};
			
			\draw (2.5,0.5) node[draw]{$1$};
			\draw (2.5,1.5) node[draw]{$0$};
			\draw (2.5,2.5) node[draw]{$0$};
			\draw (2.5,3.5) node[draw]{$0$};
			\draw (2.5,4.5) node[draw]{$0$};
			
			\draw (3.5,0.5) node[draw]{$1$};
			\draw (3.5,1.5) node[draw]{$0$};
			\draw (3.5,2.5) node[draw]{$0$};
			\draw (3.5,3.5) node[draw]{$0$};
			\draw (3.5,4.5) node[draw]{$1$};
			
			\draw (4.5,0.5) node[draw]{$1$};
			\draw (4.5,1.5) node[draw]{$1$};
			\draw (4.5,2.5) node[draw]{$1$};
			\draw (4.5,3.5) node[draw]{$1$};
			\draw (4.5,4.5) node[draw]{$0$};
		\end{tikzpicture}
    \caption{Proposition 7}
	\end{figure}

Even more restrictive operations share these peculiarities.
 
	\begin{prop}\label{pereodic7}
		An insertion  of a zero column to the left of the Bouton's matrix can change the  valuation.
For example:
    \[
		\textcolor{darkblue}{NN} \ (6, 9, 10, 11, 59) \to \textcolor{darkblue}{PN} \ (12, 18, 20, 22, 22)
		\]
	\end{prop}

	\begin{figure}[H]
		\centering
		\begin{tikzpicture}
			\draw[step=1cm,gray,very thin] (0,0) grid (5, 4);
			\draw (0.5,4.5) node[draw, red]{$NN$};
			
			\draw (0.5,0.5) node[draw]{$1$};
			\draw (0.5,1.5) node[draw]{$0$};
			\draw (0.5,2.5) node[draw]{$1$};
			\draw (0.5,3.5) node[draw]{$0$};
			
			\draw (1.5,0.5) node[draw]{$1$};
			\draw (1.5,1.5) node[draw]{$1$};
			\draw (1.5,2.5) node[draw]{$0$};
			\draw (1.5,3.5) node[draw]{$1$};
			
			%\draw (2.5,0.5) node[draw]{$0$};
			\draw (2.5,0.5) node[draw]{$0$};
			\draw (2.5,1.5) node[draw]{$0$};
			\draw (2.5,2.5) node[draw]{$0$};
			\draw (2.5,3.5) node[draw]{$1$};
			
			%\draw (3.5,0.5) node[draw]{$1$};
			\draw (3.5,0.5) node[draw]{$1$};
			\draw (3.5,1.5) node[draw]{$1$};
			\draw (3.5,2.5) node[draw]{$1$};
			\draw (3.5,3.5) node[draw]{$0$};

		\end{tikzpicture}
		\hskip1cm
		\begin{tikzpicture}
			\draw[step=1cm,gray,very thin] (0,0) grid (5, 4);
			\draw (0.5,4.5) node[draw, red]{$PN$};
			
			%\draw (0.5,0.5) node[draw]{$0$};
			\draw (0.5,0.5) node[draw]{$0$};
			\draw (0.5,1.5) node[draw]{$0$};
			\draw (0.5,2.5) node[draw]{$0$};
			\draw (0.5,3.5) node[draw]{$0$};
			
			%\draw (1.5,0.5) node[draw]{$1$};
			\draw (1.5,0.5) node[draw]{$1$};
			\draw (1.5,1.5) node[draw]{$0$};
			\draw (1.5,2.5) node[draw]{$1$};
			\draw (1.5,3.5) node[draw]{$0$};
			
			%\draw (2.5,0.5) node[draw]{$1$};
			\draw (2.5,0.5) node[draw]{$1$};
			\draw (2.5,1.5) node[draw]{$1$};
			\draw (2.5,2.5) node[draw]{$0$};
			\draw (2.5,3.5) node[draw]{$1$};
			
			%\draw (3.5,0.5) node[draw]{$0$};
			\draw (3.5,0.5) node[draw]{$0$};
			\draw (3.5,1.5) node[draw]{$0$};
			\draw (3.5,2.5) node[draw]{$0$};
			\draw (3.5,3.5) node[draw]{$1$};
			
			%\draw (4.5,0.5) node[draw]{$1$};
			\draw (4.5,0.5) node[draw]{$1$};
			\draw (4.5,1.5) node[draw]{$1$};
			\draw (4.5,2.5) node[draw]{$1$};
			\draw (4.5,3.5) node[draw]{$0$};
			
		\end{tikzpicture}
  \caption{Proposition 8}
	\end{figure}
	
	\begin{prop}\label{pereodic8}
		There is a pair of $PN$ positions, which differs by 
  	 insertion  of a zero column to the left of the reduced matrix  $\hx$. For example:
  
		\[
		  \{ (20, 33, 36, 37, 37), (40, 66, 72, 74, 74)\}
		\]
	\end{prop}
	
	\begin{figure}[!h]
		\centering
		\begin{tikzpicture}
			\draw[step=1cm,gray,very thin] (0,0) grid (7, 5);
			\draw (0.5,5.5) node[draw, red]{$PN$};
			
			\draw (0.5,0.5) node[draw]{$1$};
			\draw (0.5,1.5) node[draw]{$1$};
			\draw (0.5,2.5) node[draw]{$0$};
			\draw (0.5,3.5) node[draw]{$1$};
			\draw (0.5,4.5) node[draw]{$0$};
			
			\draw (1.5,0.5) node[draw]{$0$};
			\draw (1.5,1.5) node[draw]{$0$};
			\draw (1.5,2.5) node[draw]{$0$};
			\draw (1.5,3.5) node[draw]{$0$};
			\draw (1.5,4.5) node[draw]{$0$};
			
			\draw (2.5,0.5) node[draw]{$1$};
			\draw (2.5,1.5) node[draw]{$1$};
			\draw (2.5,2.5) node[draw]{$1$};
			\draw (2.5,3.5) node[draw]{$0$};
			\draw (2.5,4.5) node[draw]{$1$};
			
			\draw (3.5,0.5) node[draw]{$0$};
			\draw (3.5,1.5) node[draw]{$0$};
			\draw (3.5,2.5) node[draw]{$0$};
			\draw (3.5,3.5) node[draw]{$0$};
			\draw (3.5,4.5) node[draw]{$0$};
			
			\draw (4.5,0.5) node[draw]{$0$};
			\draw (4.5,1.5) node[draw]{$0$};
			\draw (4.5,2.5) node[draw]{$0$};
			\draw (4.5,3.5) node[draw]{$0$};
			\draw (4.5,4.5) node[draw]{$1$};
			
			\draw (5.5,0.5) node[draw]{$1$};
			\draw (5.5,1.5) node[draw]{$1$};
			\draw (5.5,2.5) node[draw]{$1$};
			\draw (5.5,3.5) node[draw]{$1$};
			\draw (5.5,4.5) node[draw]{$0$};
			
		\end{tikzpicture}
		\hskip1cm
		\begin{tikzpicture}
			\draw[step=1cm,gray,very thin] (0,0) grid (7, 5);
			\draw (0.5,5.5) node[draw, red]{$PN$};
			
			\draw (0.5,0.5) node[draw]{$0$};
			\draw (0.5,1.5) node[draw]{$0$};
			\draw (0.5,2.5) node[draw]{$0$};
			\draw (0.5,3.5) node[draw]{$0$};
			\draw (0.5,4.5) node[draw]{$0$};
			
			\draw (1.5,0.5) node[draw]{$1$};
			\draw (1.5,1.5) node[draw]{$1$};
			\draw (1.5,2.5) node[draw]{$0$};
			\draw (1.5,3.5) node[draw]{$1$};
			\draw (1.5,4.5) node[draw]{$0$};
			
			\draw (2.5,0.5) node[draw]{$0$};
			\draw (2.5,1.5) node[draw]{$0$};
			\draw (2.5,2.5) node[draw]{$0$};
			\draw (2.5,3.5) node[draw]{$0$};
			\draw (2.5,4.5) node[draw]{$0$};
			
			\draw (3.5,0.5) node[draw]{$1$};
			\draw (3.5,1.5) node[draw]{$1$};
			\draw (3.5,2.5) node[draw]{$1$};
			\draw (3.5,3.5) node[draw]{$0$};
			\draw (3.5,4.5) node[draw]{$1$};
			
			\draw (4.5,0.5) node[draw]{$0$};
			\draw (4.5,1.5) node[draw]{$0$};
			\draw (4.5,2.5) node[draw]{$0$};
			\draw (4.5,3.5) node[draw]{$0$};
			\draw (4.5,4.5) node[draw]{$0$};
			
			\draw (5.5,0.5) node[draw]{$0$};
			\draw (5.5,1.5) node[draw]{$0$};
			\draw (5.5,2.5) node[draw]{$0$};
			\draw (5.5,3.5) node[draw]{$0$};
			\draw (5.5,4.5) node[draw]{$1$};
			
			\draw (6.5,0.5) node[draw]{$1$};
			\draw (6.5,1.5) node[draw]{$1$};
			\draw (6.5,2.5) node[draw]{$1$};
			\draw (6.5,3.5) node[draw]{$1$};
			\draw (6.5,4.5) node[draw]{$0$};
			
		\end{tikzpicture}
        \caption{Proposition 9}
	\end{figure}

    \section{Exceptional positions}

For better understanding of difference between the exact $(5,2)$ nim and the Moore's $(4,2)$ nim we introduce another classification of positions in good, bad and exceptional.
%It would be useful to reduce the number of the positions and study the small class of them, thus, we defined exceptional positions.

The definition is by the backward induction. Assume that for a position $x$ all  moves go to  positions with known result of classification (either good, or bad, or exceptional).
	
	%Storing three sets:
	
	%\begin{enumerate}
	%	\item Good
	%	\item Bad
	%	\item Exceptional
	%\end{enumerate}
	
	Then the classification of $x$ is determined by the following rules.
	\begin{itemize}
		\item If $x$ is PP then $x$ is good.
		\item If $x$ is NP then $x$ is bad.
		\item If $x$ is NN and there is a move from $x$ to either PP or NP, then $x$ is good.
		\item if $x$ is PN and there is a move from $x$ to  either PP or NP, then $x$ is bad.
		\item If $x$ is PN and there is no move to PP or NP, then $x$ is good.
		\item If $x$ is NN and
  there is no move to PP or NP, then $x$ is bad.
		\item If an N-position $x$ in the exact (5,2) nim is bad, but there is a move to a good P-position, then $x$  is regular. Otherwise, (there is no move to a good P) it is exceptional.
		\item If a P-position in the exact (5,2) nim is bad, but all moves go to good N-positions, it is regular. Otherwise, (there is a move to a bad N) it is exceptional.  
	\end{itemize}

 The file of  exceptional positions contained in the computed list of positions can be downloaded via link \cite{exceptionalFILE}.
 %File with exceptional among our file can be downloaded via link. \cite{exceptionalFILE}
Moves between exceptional positions form a digraph called the exceptional graph.

We have found the following observations on the computed exceptional positions.
 
	\begin{conj}\label{exceptional1}
		For all exceptional N-positions, the outgoing degree in the exceptional graph is multiple of $3$.   
	\end{conj}
	
	\begin{conj}\label{exceptional2}
		All exceptional N-positions are NP. 

		Because of that, there is an alternation on exceptional positions.
		All terminals on the exceptional graph are P-positions, there is a move to them
		from N-positions, which is $\PM{}$, to which there is a move from 
		P-positions and so on.

		The positions $\NM{}$, from which there is no move to $\PM{}$ - we call \textbf{deadenders}. 
		
		Besides, to that alternating chains we should add NP-positions, from
		which there is a move to deadenders. But such positions are not exceptional. 
	\end{conj}
	
\begin{figure}
        \includegraphics[scale=0.85]{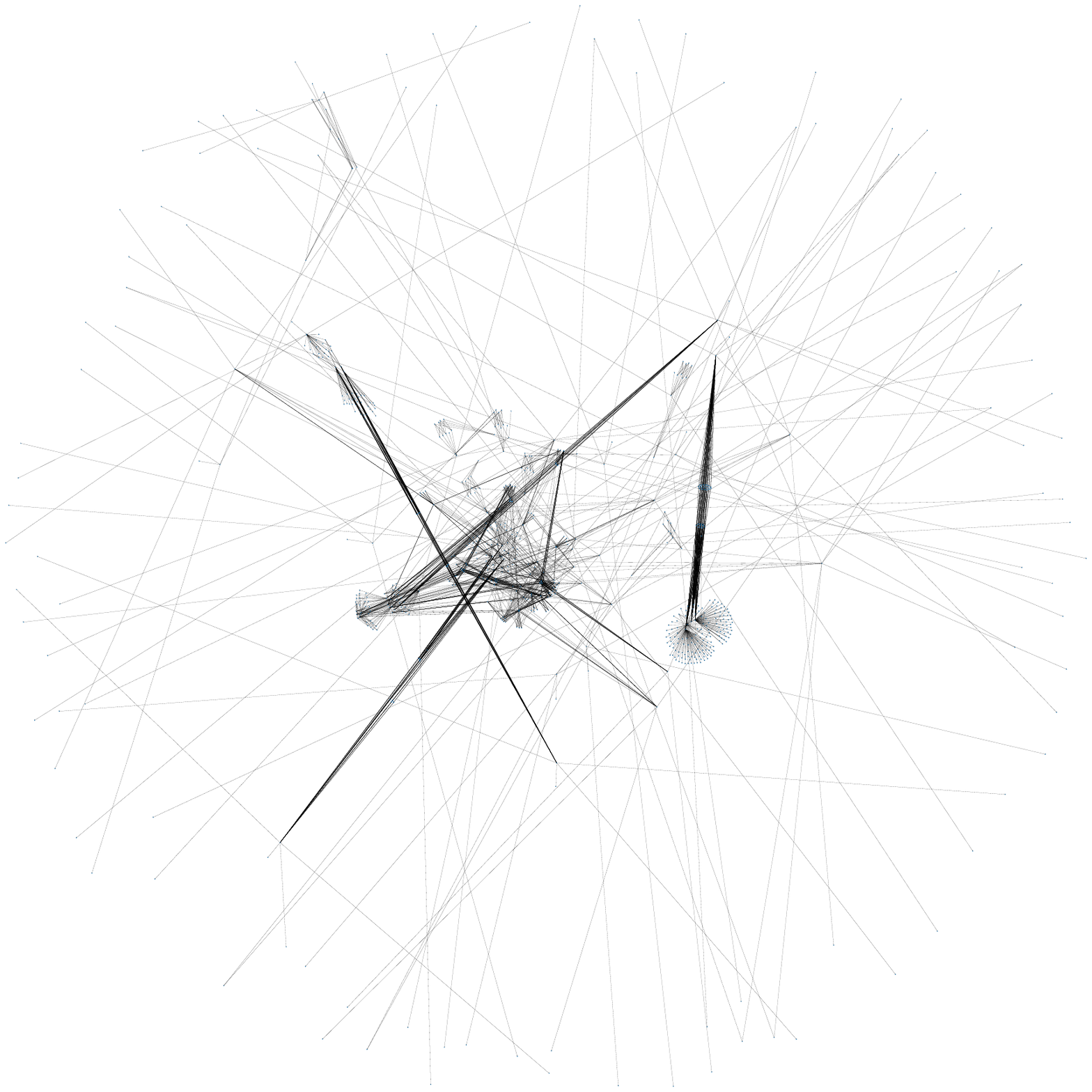}
        \caption{Graph of the game on exceptional positions without isolated vertexes}\label{pic3}
    \end{figure}

	On figure \ref{pic3}, p.~\pageref{pic3}, you can see the graph of the game on exceptional positions without isolated vertexes. Most of the moves are concentrated near few positions. 
	Also, in that graph there are a few weak connectivity components with big cardinality, but the size of such components decreases very fast.

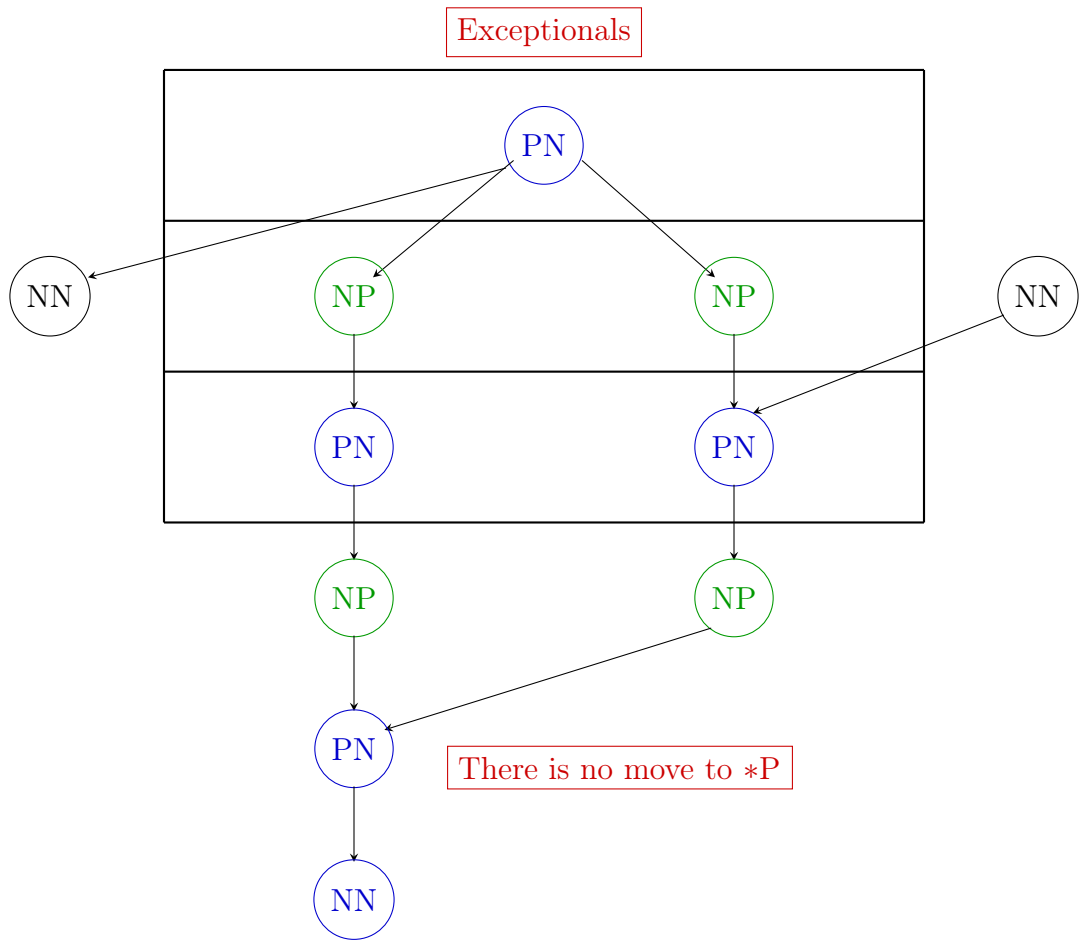
\begin{figure}[H]
        \centering
        \begin{tikzpicture}
		
		\draw (5,10.5) node[draw, darkred]{Exceptionals};
		\draw (5, 9) node[draw, circle, darkblue]{PN};
		
		\draw[thick,-] (0,10) -- (10,10);
		\draw[thick,-] (0,8) -- (10,8);
		\draw[thick,-] (0,6) -- (10,6);
		\draw[thick,-] (0,4) -- (10,4);
		%\draw[thick,-] (0,2) -- (10,2);
		%\draw[thick,-] (0,0) -- (10,0);
		
		\draw [-stealth](4.6, 8.8) -- (2.75, 7.25);
		\draw [-stealth](5.5, 8.8) -- (7.25, 7.25);
		\draw (-1.5, 7) node[draw, circle, black]{NN};
		\draw [-stealth](4.5, 8.7) -- (-1, 7.25);
		
		\draw (7.5, 7) node[draw, circle, darkgreen]{NP};
		\draw (2.5, 7) node[draw, circle, darkgreen]{NP};
		
		\draw (11.5, 7) node[draw, circle, black]{NN};
		\draw [-stealth](11.05, 6.75) -- (7.75, 5.45);
		
		\draw [-stealth](2.5, 6.5) -- (2.5, 5.5);
		\draw [-stealth](7.5, 6.5) -- (7.5, 5.5);
		
		\draw (2.5, 5) node[draw, circle, darkblue]{PN};
		\draw (7.5, 5) node[draw, circle, darkblue]{PN};
		
		\draw [-stealth](2.5, 4.5) -- (2.5, 3.5);
		\draw [-stealth](7.5, 4.5) -- (7.5, 3.5);
		
		\draw (7.5, 3) node[draw, circle, darkgreen]{NP};
		\draw (2.5, 3) node[draw, circle, darkgreen]{NP}; 
		
		\draw [-stealth](2.5, 2.5) -- (2.5, 1.5);
		\draw [-stealth](7.2, 2.6) -- (2.9, 1.25);
		
		\draw (2.5, 1) node[draw, circle, darkblue]{PN};
		
		\draw (6, 0.75) node[draw, darkred]{There is no move to $\PM{}$};
		
		\draw [-stealth](2.5, 0.5) -- (2.5, -0.5);
		
		\draw[thick,-] (0,10) -- (0,4);
		\draw[thick,-] (10,10) -- (10,4);
		
		\draw (2.5, -1) node[draw, circle, darkblue]{NN};

	\end{tikzpicture}	
	
        \caption{Exceptional positions}
        \label{except}
    \end{figure}

    Exceptional positions can be seen on the graph \ref{except}

	There might be more chains, it might go to different terminals and grow up infinitely, it's just a scheme.

	\section{Remoteness function}
	
	\subsection{Definition}

    The remoteness function can be defined for each game in such way.
	
	It will be evaluated recursively. 
	
	Let $s = 0$. For all terminals $\R_{n, k} = s$. For all positions, from which there is a move to terminals let $\R_{n, k} = s + 1$. 
	Then let's throw away the positions for which we already evaluated function and increase $s \to s + 2$. 
	Then repeat that on a smaller graph of the game.

Let    $ \R_{5, 2}(x)$ the value of the remoteness function of a position $x$ in the exact $(5, 2)$ nim. File with the values of this function can be downloaded via link~\cite{SmithFILE}

	The relation between the Moore's nim$(4, \leq 2)$ and the exact $(5, 2)$ nim
	can be illustrated in terms of the remoteness function. For most of the positions $x$ we have $\R_{4, 2}(\hx) = \R_{5, 2}(x)$

	Obviously, for exceptional positions they differ, moreover not always by one.
	Starting at some point difference grows.

%	Denote m-critical positions minimal in coordinate order positions with value of remoteness function equals m.

	%\begin{conj}\label{crit1}
	%	There is no critical exceptional positions.
	%\end{conj}

	\section{Properties of $NP$ positions}

	Let's recall that $i$-th column of the reduced matrix is balanced for a position $x$, if $m(x)_i\in\{0,3\}$.
	Thus, Moore's criterion claims that $x$ is $\PM{}$-position, if and only if all columns are balanced.

	\begin{defi}
 Let $\xi_w(x)$ be the integer such that $k$-th bit of its  binary representation  is~1, if $M(x)_k = w$.
		In other words, 
			\[
			\xi_w(x) = \sum_k [M(x)_{k}=w]\cdot2^k.
			\]
	\end{defi}

	From the definition it is obvious that for $w'\ne w''$ ones in binary representations 
	$\xi_{w'}(x)$ and $\xi_{w''}(x)$ stands for different bits. 

For positions in the computed list the following holds.
 
	\begin{conj}%[Criterion of  NP positions]
 \label{th:NP}
		The position $x$ is \textup{NP}, if and only if $M(x)\in\{0,3\}^*$ and $\xi_3(x)> x_5$. 
	\end{conj}

	We conjecture that it always holds.
	 
\noindent  \textbf{Conjecture.} PP positions are positions so $M(x)\in\{0,3\}^*$ and $\xi_3(x)\leq x_5$.
	
	Just now, we can't prove this claim due to the lack of plausible conjectures about other classes of positions. But, in our opinion, the conjecture has a strong support since  there are no counterexamples for it in the computed list.

	\section{Proofs}

    To extend this report by several proofs. Some of them are known, some of them are not, and an interested reader can read them here.

 %    \subsection{The criterion of the P-positions in the nim}

 %    Let's recall the theorem

 %    \begin{theorem}[Criterion of P-positions]
	% 	The position $x$ in the nim is \textup{P}-position, if and only if
	% 	\[
	% 	\bigoplus\limits_{i = 1}^{n} x_{ij} = 0\quad\text{for each $j$.}
	% 	\]
 %    \end{theorem}

 %    \begin{proof}[Proof]
	% 	Induction.
		
	% 	Base: Terminal positions --- all piles are empty --- done.
		
	% 	From P-positions all moves lead to N-positions. Indeed, 
	% 	if the sum is zero, let's take from the $i$-th pile $t$ stones. 
	% 	Then the sum will be equal to $x_i \oplus (x_i - t) \neq 0$.
		
	% 	From N-positions we can make a move to the P-position (at least one, 
	% 	but in general maybe there are several such moves). Indeed, 
	% 	it is always possible to make a zero-sum from a non-zero sum by 
	% 	taking $x_k$ such that the highest bit of the sum is not zero 
	% 	(let's call it $j$). Then by reducing the number of stones 
	% 	to $\bigoplus\limits_{i=1}^{n}x_i\oplus x_k$.

	% 	If the highest bit of the total sum is not zero, then the same bit 
	% 	must be equal to one for one (or any other odd number) of piles. At the
	% 	same time, reducing the number of stones in the pile to 
	% 	$\bigoplus\limits_{i = 1}^{n} x_i\oplus x_k$ is exactly a valid move,
	% 	since we knowingly add two numbers whose $j$-th bit is equal to one, 
	% 	therefore, even if all the previous bits become equal to one, 
	% 	the number of stones in the pile will decrease according to the 
	% 	properties of binary numbers. Thus, such a heap must exist.

	% \end{proof}

    \subsection{The Moore's criterion}

    Let's repeat the criterion  of P-positions for the Moore's nim$(n, \leq k)$: 
		$x$ is a P-position if and only if 
		$			M(x)_j\equiv 0\pmod{k+1}
		$.

    \begin{proof}[Proof]
		Induction on the sum of coordinates. 
		
		The induction base: if the sum is 0, then definitely $M(x) = 0$.

  The induction step. We assume that for the sum of coordinates less than $S$ the criterion holds. Consider a position $x$ with the sum of coordinates equals $S$. Note that every move decreases the sum of coordinates. Thus by induction hypothesis the criterion holds for all positions reachable from $x$ in one move.
		\begin{itemize}
			\item[$\Longrightarrow$] Suppose  that $M(x)_j\equiv 0\pmod{k+1}$ for all $j$. 
   %if $x$ is a P-position then  all moves from $x$ go to N-positions.
			
			%In P-positions, sum $\pmod{k + 1}$ equals zero. 
   In one move $x\to x'$ a player takes stones from at most $k$ piles. So, in each column of Bouton's matrix, at most $k$ bits are changed after a move, and there exists a column $j$, in which at least one bit is changed. Thus, after a move, $M(x')_j\not\equiv 0\pmod {k+1}$. Due to the induction hypothesis, it means that $x'$ is an N-position. Therefore $x$ is a P-position.
%   Because for each move we can change the value
	%		in less than $k + 1$ piles, in order to get after the move non-zero value $\pmod{k + 1}$, we need for 
		%	each bit, for which we took, change that corresponding bit in some other pile.

			%Now let's look at the first from the left-hand side changed bit(we use the little-endian notation). Since we took 
			%some bit in this digit, it means that in one of the other piles we have to replace the zero bit with one in order to keep the balance.
			
%			But we can only take --- not add --- some stones, so in order to replace some bit from zero to one, we need to make reverse change in some lefter bit
			%But this contradicts with the fact that we were looking at the first from the left changed bit.

			\item[$\Longleftarrow$]
   Suppose now  that $M(x)_j\not \equiv 0\pmod{k+1}$ for some $j$.
   To prove that $x$ is an N-position, we should indicate a move  to a P-position $x'$. By the induction hypothesis, it is equivalent to the existence of a move $x\to x'$ such that $M(x')_i \equiv 0\pmod{k+1}$ for all $i$.
			
			Let $j$ be the highest unbalanced column. Our goal is to decrease some bits in unbalanced columns starting from $j$th  to make $M(x')_i$ zero modulo $k + 1$.

We are going through the columns one by  one. Let $t$ be the  number of rows bearing the  changed bits before processing an unbalanced column $i$.  We call these rows changed. If, in the current column $i$, it is possible to make changes in at most $ t$ changed rows to make the value $M(x')_i$ equal to zero, we do so. This can happen if $t > 0$, which implies that in higher columns  we already changed something. Thus, both changes $0\to1$ and $1\to 0$ are possible.

Suppose that,  among $t$ changed rows, there are $u$ ones and the sum in the column $i$ modulo ${k + 1} $ is $ a$.
If $a\leq u$ the above operation is possible. So, in the sequel, we assume that $a>u$.

	%		Now we need to find the leftest unbalanced column and let's decrease each column starting from it to zero $\pmod{k + 1}$. Also, we need to store
	%		the number of columns we have changed --- call it $t$. If for the current column it is enough to change $\leqslant t$ bits to make it equal zero, let's do so.

		%	That case can happen if $t > 0$, which means that in lefter columns we already have changed something and that's why we can change these $\leqslant t$ bits as we want.

			If we need to change more than $t$ bits to balance the column $i$, we need to add another $t'$ rows to the set of changed rows such that $t + t' \leqslant k$ and only change $1\to 0$ is possible in these additional rows. 
			After  that, we have  $t + t'$ changed rows.

			If $t + a - u \leqslant k$, then we can make $M(x')_i$ zero modulo $k + 1$ by replacing in $t$
			changed rows all ones to zeros and additional $a - u$ ones to zeros in other rows. Thus $t' = a-u$ and $t+t' = t+a-u\leq k$, as required.
			
			Otherwise, $t + a - u > k$ and we need to add ones. Among $t$ rows changed before this step we change $k + 1 - a$ zeros to ones (the total amount of zeros in these rows is $t - u$).
			It is possible if $k + 1 - a \leqslant t - u $, equivalently, $  k + 1 \leqslant t - u + a$ and it is equivalent to the assumption for this case. After that, the column $i$ becomes balanced and $t'=0$. 
	   In the result, we are able to balance all columns making the above operations and define a move to a P-position $x'$. Therefore $x$ is an N-position.
		\end{itemize}
\end{proof}	
			%The sums in all columns were $\pmod{k + 1}$, so the max possible $t$ after all algorithm is equal, $k$ and because of that we always can make a move to the P-position.
			
			\begin{rem}
				The classic nim is one of the cases of the Moore's nim, so this proof also works for it. 
			\end{rem}

	\subsection{Special case of the exact (5, 2) nim}

 %TODO:: Insert the reference to Hypergraph paper (or other with proof of exact 2n, n)

    In fact, this is just the exact (4, 2) nim. We reproduce the criterion of P-positions from~\cite{BOROS20181}.
	
	\begin{theorem}[]\label{P4-2}
				P-positions $(0,a,b,c,d)$ are the positions
		\[(0,a,a,a,b).\]
	\end{theorem}
	
	\begin{proof}[Proof]
		Induction.

		Basis: Terminals indeed looks like that. 
		
		\begin{itemize}
			\item[$\Longrightarrow$] Let's show that from any N-position, 
			there is a move to P-position.
			\begin{enumerate}
				\item $x_2 < x_3$: the move is to decrease $x_3, x_4$ to $x_2$. 
				\item $x_2 = x_3$: decrease $x_4, x_5$ to $x_2$.
			\end{enumerate}
			
			\item[$\Longleftarrow$]  From  \[
			(0, a, a, a, b)
			\]
			all moves goes to either positions for which there is no 3 equal piles, 
			or to the positions with different central\footnote{not first and not last} piles.
			
			Indeed, we have two options:
			
			\begin{enumerate}
				\item Take from two piles $a$ stones. Such way we can go to 
				one of two cases:
				\begin{enumerate}
					\item $(0, t, t, a, b), \ 0 < t < a < b$
					\item $(0, s, t, a, b), \ 0 < s < t < a < b$
				\end{enumerate}
				
				In both we don't have equal central piles.
				
				\item Take from one pile with $a$ stones and from pile with
				 $b$ stones. Such way we can go to:
				\begin{enumerate}
					\item $(0, t, a, a, a), \ 0 < t < a$ (Different number of stones in central piles).
					\item $(0, t, t, a, a), \ 0<  t < a$ (Less than three equal piles).
					\item $(0, t, a, a, s), \ 0 < t < a < s$ (Less than three equal piles).
				\end{enumerate}
			\end{enumerate}
			
			Thus, from P-position, it is possible to go only to N-positions.
		\end{itemize}

	\end{proof}
	
	\subsection{Move $\PM{} \to \PM{}$ in the exact (5, 2) nim}

        There is an interesting fact about the moves from $\PM{}$ to $\PM{}$. We prove that such moves don't exist with the following argument.

        Obviously there are no moves $\PM{} \to \PM{}$ without changing the greatest pile,
		otherwise there exists f move from P-position to P-position in the Moore's nim $(4, \leq 2)$, a contradiction.
 
	\begin{lemma}[Move $\PM{} \to \PM{}$]

		There is no move in the exact $(5, 2)$ nim from $\PM{}$ to $\PM{}$ with changing the leader (pile with maximum stones in position).
	\end{lemma}
	%%    \proofname. 
	\begin{proof}[Proof]
		By contradiction. 
		
		Suppose $x\to y$ is a move from $\PM{}$ to $\PM{}$.

		We can change the leader only if we take some stones from previous leader.

		Denote a row in the reduced matrix with the previous leader as $o$, a row with a new leader as $n$.
		
		Since $y_n>y_o$, there exists $k$ such that $y_{nk}=1$, $y_{ok} =0$ and $y_{nj} = y_{oj}$ for $j>k$.

		\begin{enumerate}
			\item Suppose  $M(y)_{k} = 3$. 
			First, suppose that $x_{nk} = 1$.
			%$ \wedge x_{ok} = 0$.
		The new leader is indicated with thick parallel lines on the figure. 
   %Let's look at the column, because of which we have changed the leader:
	
 In the $k$-th column of the reduced matrix  $\hx$ there are exactly three ones (the column is balanced and contains a one by assumption of this case). Thus $x_{ak}=0$ for some $a\ne o$. Then $y_{ak} =1$, since $M(y)_{k} = 3$ and $y_{ok}=0$. The change $x_{ak} = 0 \to y_{ak} = 1$ during the move is possible only if in some column $t > k$ the change $x_{at} = 1 \to 0 = y_{at}$ occurs. In particular, it means that the move takes stones from exactly rows $o$, $a$.

			\medskip
			
			\begin{tikzpicture}
				\draw[step=1cm,gray,very thin] (0,0) grid (5,5);
				\draw[line width=0.5mm,-] (0,1) -- (5,1);
				\draw[line width=0.5mm,-] (0,2) -- (5,2);
				\draw (0.5,5.5) node[draw, red]{$y$};
				\draw (1.5,1.5) node[draw]{$1$};
				\draw (1.5,0.5) node[draw]{$0$};
				\draw (1.5,2.5) node[draw]{$1$};
				\draw (1.5,3.5) node[draw]{$1$};
				\draw (1.5,4.5) node[draw]{$1$};
				\draw (1.5,5.5) node[blue]{$k$};
                \draw (3.5,1.5) node[draw]{$u$};
				\draw (3.5,0.5) node[draw]{$u$};
				\draw (3.5,2.5) node[circle, draw]{$0$};
				\draw (3.5,3.5) node[draw]{$1$};
				\draw (3.5,4.5) node[draw]{$1$};
				\draw (1.5,5.5) node[blue]{$k$};
				\draw (3.5,5.5) node[blue]{$t$};

                \draw (-0.2,2.5) node[blue]{$a$};
                \draw (-0.2,1.5) node[blue]{$n$};
                \draw (-0.2,0.5) node[blue]{$o$};
			\end{tikzpicture}
			\hskip1cm
			\begin{tikzpicture}
				\draw[step=1cm,gray,very thin] (0,0) grid (5,5);
				\draw[line width=0.5mm,-] (0,1) -- (5,1);
				\draw[line width=0.5mm,-] (0,2) -- (5,2);
				\draw (0.5,5.5) node[draw, red]{$x$};
				\draw (1.5,1.5) node[draw]{$1$};
				\draw (1.5,0.5) node[draw]{$0$};
				\draw (1.5,2.5) node[draw]{$0$};
				\draw (1.5,3.5) node[draw]{$1$};
				\draw (1.5,4.5) node[draw]{$1$};
				\draw (3.5,1.5) node[draw]{$x_{nt}$};
				\draw (3.5,0.5) node[draw]{$x_{ot}$};
				\draw (3.5,2.5) node[circle, draw]{$1$};
				\draw (3.5,3.5) node[draw]{$1$};
				\draw (3.5,4.5) node[draw]{$1$};
				\draw (1.5,5.5) node[blue]{$k$};
				\draw (3.5,5.5) node[blue]{$t$};

                \draw (-0.2,2.5) node[blue]{$a$};
                \draw (-0.2,1.5) node[blue]{$n$};
                \draw (-0.2,0.5) node[blue]{$o$};
			\end{tikzpicture}
	
Note that in the remaining rows that differ from $a,o,n$ the bits are not changed because the rows $o$ and $a$ are necessarily changed. Therefore they contain ones in $k$-th column as well as in $t$-th column, since $M(x)_t > 0$ (thus $M(x)_t = 3)$ and $M(y)_k =3$ by the assumption. 

Note that $y_{nt} = y_{ot} = 1$, since $t$-th column of $y$ is balanced and $t>k$. But $x_{nt} = y_{nt} =1 $, since the move takes stones from $a$ and $o$. Thus, the $t$-th column in $\hat x$ is unbalanced, it contains four ones. A contradiction.

			\item 
			Assume  $M(y)_{k} = 3$ and $x_{nk} = 0$. Again, the change $x_{nk} = 0 \to 1 = y_{nk}$ implies that 
			for some $t>k$ the change $x_{nt} = 1 \to 0 = y_{nt}$ occurs. Thus, the move takes stones from $o$, $n$ and bits in the rest of rows are not changed. Since $x_{nt}=1$ and $t$-th column in $\hat x$ is balanced, there is a row $a$, $a\ne n$, $a\ne o$, such that $x_{at} =0 = y_{at}$. We come to a contradiction with the assumption $M(y)_{k} = 3$.

			\medskip
			
			\begin{tikzpicture}
				\draw[step=1cm,gray,very thin] (0,0) grid (5,5);
				\draw[line width=0.5mm,-] (0,1) -- (5,1);
				\draw[line width=0.5mm,-] (0,2) -- (5,2);
				\draw (0.5,5.5) node[draw, red]{$y$};
				\draw (1.5,1.5) node[draw]{$1$};
				\draw (1.5,0.5) node[draw]{$0$};
				\draw (1.5,2.5) node[draw]{$1$};
				\draw (1.5,3.5) node[draw]{$1$};
				\draw (1.5,4.5) node[draw]{$1$};
				\draw (1.5,5.5) node[blue]{$k$};
				\draw (3.5,5.5) node[blue]{$t$};
				\draw (3.5,1.5) node[circle, draw]{$0$};
				\draw (3.5,0.5) node[draw]{$0$};
				\draw (3.5,2.5) node[draw]{$1$};
				\draw (3.5,3.5) node[draw]{$0$};
				\draw (3.5,4.5) node[draw]{$1$};

                \draw (-0.2,1.5) node[blue]{$n$};
                \draw (-0.2,0.5) node[blue]{$o$};
                \draw (-0.2,3.5) node[blue]{$a$};

			\end{tikzpicture}
			\hskip1cm
			\begin{tikzpicture}
				\draw[step=1cm,gray,very thin] (0,0) grid (5,5);
				\draw[line width=0.5mm,-] (0,1) -- (5,1);
				\draw[line width=0.5mm,-] (0,2) -- (5,2);
				\draw (0.5,5.5) node[draw, red]{$x$};
				\draw (1.5,1.5) node[draw]{$0$};
				\draw (1.5,0.5) node[draw]{$0$};
				\draw (1.5,2.5) node[draw]{$1$};
				\draw (1.5,3.5) node[draw]{$1$};
				\draw (1.5,4.5) node[draw]{$1$};
				\draw (3.5,1.5) node[circle, draw]{$1$};
				\draw (3.5,0.5) node[draw]{$1$};
				\draw (3.5,2.5) node[draw]{$1$};
				\draw (3.5,3.5) node[draw]{$0$};
				\draw (3.5,4.5) node[draw]{$1$};
				\draw (1.5,5.5) node[blue]{$k$};
				\draw (3.5,5.5) node[blue]{$t$};

                \draw (-0.2,1.5) node[blue]{$n$};
                \draw (-0.2,0.5) node[blue]{$o$};
                \draw (-0.2,3.5) node[blue]{$a$};
			\end{tikzpicture}
			
			\item 
			Assume  $M(y)_{k} = 0$.
   Then $M(x)_k=0$ also, since otherwise at least two bits in rows that differ from $n$, $o$ are changed from 1 to 0. It is impossible because the move does take stones from $o$. In particular, $x_{nk}=0$. Note that $y_{nk}=1 $ by construction. It implies that the move takes stones from rows $n$, $o$. 

   Since $x_{nt} = 0$ and $y_{nt}=1$, 
  		 for some $t > k$ we have $x_{nt} = 1$, $y_{nt}=0$. 
			
			\begin{tikzpicture}
				\draw[step=1cm,gray,very thin] (0,0) grid (5,5);
				\draw[line width=0.5mm,-] (0,1) -- (5,1);
				\draw[line width=0.5mm,-] (0,2) -- (5,2);
				\draw (0.5,5.5) node[draw, red]{$y$};
				\draw (1.5,1.5) node[draw]{$1$};
				\draw (1.5,0.5) node[draw]{$0$};
				\draw (1.5,2.5) node[draw]{$0$};
				\draw (1.5,3.5) node[draw]{$0$};
				\draw (1.5,4.5) node[draw]{$0$};
				\draw (1.5,5.5) node[blue]{$k$};
				\draw (3.5,5.5) node[blue]{$t$};
				\draw (3.5,1.5) node[circle, draw]{$0$};
				\draw (3.5,0.5) node[draw]{$0$};
				\draw (3.5,2.5) node[draw]{$1$};
				\draw (3.5,3.5) node[draw]{$0$};
				\draw (3.5,4.5) node[draw]{$1$};

                \draw (-0.2,1.5) node[blue]{$n$};
                \draw (-0.2,0.5) node[blue]{$o$};
                \draw (-0.2,2.5) node[blue]{$c$};
                \draw (-0.2,3.5) node[blue]{$b$};
                \draw (-0.2,4.5) node[blue]{$a$};
			\end{tikzpicture}
			\hskip1cm
			\begin{tikzpicture}
				\draw[step=1cm,gray,very thin] (0,0) grid (5,5);
				\draw[line width=0.5mm,-] (0,1) -- (5,1);
				\draw[line width=0.5mm,-] (0,2) -- (5,2);
				\draw (0.5,5.5) node[draw, red]{$x$};
				\draw (1.5,1.5) node[draw]{$0$};
				\draw (1.5,0.5) node[draw]{$0$};
				\draw (1.5,2.5) node[draw]{$0$};
				\draw (1.5,3.5) node[draw]{$0$};
				\draw (1.5,4.5) node[draw]{$0$};
				\draw (3.5,1.5) node[circle, draw]{$1$};
				\draw (3.5,0.5) node[draw]{$1$};
				\draw (3.5,2.5) node[draw]{$1$};
				\draw (3.5,3.5) node[draw]{$0$};
				\draw (3.5,4.5) node[draw]{$1$};
				\draw (1.5,5.5) node[blue]{$k$};
				\draw (3.5,5.5) node[blue]{$t$};

                \draw (-0.2,1.5) node[blue]{$n$};
                \draw (-0.2,0.5) node[blue]{$o$};
                \draw (-0.2,2.5) node[blue]{$c$};
                \draw (-0.2,3.5) node[blue]{$b$};
                \draw (-0.2,4.5) node[blue]{$a$};
                
			\end{tikzpicture}

From $x_{nt} =1 $ we conclude that $M(x)_t = 3$. Thus, $x_{at} + x_{bt} + x_{ct} = 2 $. It follows from $y_{nt} = y_{ot}=0$ (by definition of $k$). But $ y_{at} + y_{bt} + y_{ct}=2$ since these bits do not change during the move. 
Therefore $M(y)_t = 2$ and we come to a contradiction with the assumption that all columns of $\hat y$ are balanced.
    \qedhere
		\end{enumerate}
	\end{proof}

 \paragraph*{Acknowledgements.} 
 This research was supported by Russian Science Foundation, 
 grant  20-11-20203, https://rscf.ru/en/project/20-11-20203/
	
\printbibliography

@article{Bouton,
author = {C.L. Bouton}, 
title = {Nim, a game with a complete mathematical theory}, 
journal = {Ann. of Math., 2-nd Ser. 3}, 
year = {1901-1902}, 
pages = {35-39}
}

@article{Moore,
author = {E. H. Moore}, 
title = {A generalization of the game called Nim}, 
journal = {Annals of Math., Second Series, 11:3}, 
year = {1910}, 
pages = {93–94}
}

@online{BaseFILE,
author = {Artem Parfenov and Mikhail Vyalyi},
title  = {List of positions with winners},
date   = {2022-11},
url    = {https://disk.yandex.ru/d/gsYUuPxp7FxV9w}
}

@online{exceptionalFILE,
author = {Artem Parfenov and Mikhail Vyalyi},
title  = {List of exceptional positions},
date   = {2023-02},
url    = {https://disk.yandex.ru/d/t3cKqGFCF7NhhQ}
}

@online{SmithFILE,
author = {Artem Parfenov and Mikhail Vyalyi},
title  = {List of Smith functions for positions},
date   = {2023-05},
url    = {https://disk.yandex.ru/d/R9_bWjl3fK6Bgw}
}

@article{BOROS20181,
title = {On the {S}prague-{G}rundy function of Exact k-Nim},
journal = {Discrete Applied Mathematics},
volume = {239},
pages = {1-14},
year = {2018},
issn = {0166-218X},
doi = {https://doi.org/10.1016/j.dam.2017.08.007},
url = {https://www.sciencedirect.com/science/article/pii/S0166218X17303931},
author = {Endre Boros and Vladimir Gurvich and Nhan Bao Ho and Kazuhisa Makino and Peter Mursic},
}

\end{document}